\documentclass[11pt,twoside]{article}

\usepackage{a4wide}
\usepackage{amsfonts}
\usepackage{amssymb}
\usepackage{amsmath}
\usepackage{graphicx}

\newcommand{\ignore}[1]{}

\def\@begintheorem#1#2{\par\bgroup{\sc #1\ #2. }\it\ignorespaces}
\def\@opargbegintheorem#1#2#3{\par\bgroup{\sc #1\ #2\ (#3). } \it\ignorespaces}
\def\@endtheorem{\egroup}
\newtheorem{theorem}{Theorem}[section]
\newtheorem{corollary}[theorem]{Corollary}
\newtheorem{lemma}[theorem]{Lemma}

\newtheorem{example}[theorem]{Example}
\newtheorem{proposition}[theorem]{Proposition}
\newtheorem{definition}[theorem]{Definition}
\newcommand{\bt}[1]{\begin{theorem}\label{#1}}
\newcommand{\bc}[1]{\begin{corollary}\label{#1}}
\newcommand{\bl}[1]{\begin{lemma}\label{#1}}
\newcommand{\be}[1]{\begin{example}\label{#1}}
\newcommand{\bp}[1]{\begin{proposition}\label{#1}}
\newcommand{\ba}[1]{\begin{algorithm}\rm\label{#1}}
\newcommand{\bd}[1]{\begin{definition}\rm\label{#1}}{\normalsize }

\newcommand{\et}{\end{theorem}}
\newcommand{\ec}{\end{corollary}}
\newcommand{\el}{\end{lemma}}
\newcommand{\ee}{\end{example}}
\newcommand{\ep}{\end{proposition}}
\newcommand{\ed}{\end{definition}}

\def\Z{\mathbb{Z}}

\begin{document}

\title{\bf Degree Sequence Optimization\\ in Bounded Treewidth}

\author{
Shmuel Onn
\thanks{\small Technion - Israel Institute of Technology. Email: onn@technion.ac.il}
\thanks{\small My manuscript has no associated data.}
}
\date{}

\maketitle

\begin{abstract}
We consider the problem of finding a subgraph of a given graph which minimizes the sum of
given functions at vertices evaluated at their subgraph degrees. While the problem
is NP-hard already when all functions are the same, we show that it can be solved
for arbitrary functions in polynomial time over graphs of bounded treewidth. Its complexity
remains widely open, in particular over complete graphs and complete bipartite graphs.

\end{abstract}

\section{Introduction}

The {\em degree sequence} of a simple graph $G=(V,E)$ on $V=[n]:=\{1,\dots,n\}$
is the vector $d(G)=(d_1(G),\dots,d_n(G))$, where $d_i(G):=|\{e\in E:i\in e\}|$
is the degree of vertex $i$ for all $i$.

\vskip.2cm
In this article we are interested in the following discrete optimization problem.

\vskip.2cm\noindent{\bf Degree Sequence Optimization.}
Given a graph $H$ on $[n]$, and for $i\in[n]$ functions
$f_i:\{0,1,\dots,n-1\}\rightarrow\Z$, find the minimum of $\sum_{i=1}^n f_i(d_i(G))$
over subgraphs $G\subseteq H$ on $[n]$.

\vskip.2cm
For a graph $G$ we let $V(G)$ and $E(G)$ be its sets of vertices and edges respectively.
Throughout, for graphs $G,H$ we use $G\subseteq H$ to indicate $V(G)=V(H)$ and $E(G)\subseteq E(H)$.

\vskip.2cm
A special case of our problem is the {\em general factor problem} \cite{Cor},
which is to decide, given $H$ and subsets $B_i\subseteq\{0,1,\dots,n-1\}$
for all $i$, if there is a $G\subseteq H$ with $d_i(G)\in B_i$ for all $i$.
Indeed, for each $i$ define $f_i(x):=0$ if $x\in B_i$ and $f_i(x):=1$
if $x\notin B_i$. Then the optimal value of our problem is zero if and
only if there is a factor. A more special case is the {\em $(l,u)$-factor problem}
introduced by Lov\'asz \cite{Lov}, where each $B_i=\{l_i,\dots,u_i\}$ is an interval.
This reduces to our problem even with convex functions, with $f_i(x):=l_i-x$ if
$0\leq x\leq l_i$, $f_i(x):=0$ if $l_i\leq x\leq u_i$, and $f_i(x):=x-u_i$ if
$u_i\leq x\leq n-1$. The {\em $b$-matching problem} is the case with $B_i=\{b_i\}$
a singleton for all $i$ and the perfect matching problem is with $b_i=1$ for all $i$.

Degree sequence optimization can be done in polynomial time when all $f_i$ are
convex \cite{AS,DO}; when all $f_i$ are the same and $H=K_n$ is complete \cite{DLMO};
and for arbitrary functions over {\em monotone} subgraphs when $H=K_{m,n}$ is complete
bipartite  \cite{KO}. For general graphs it is NP-hard even when all functions but
one are the same, as $H$ has a nonempty cubic subgraph, which is NP-complete to decide,
if and only if, for some $i$, the optimal value of degree sequence optimization
with $f_i(x):=(x-3)^2$ and $f_j(x):=x(x-3)^2$ for $j\neq i$, is equal to zero, see \cite{DO}.

\vskip.2cm
In contrast we prove here the following statement.

\bt{treewidth}
For any fixed $k$, the degree sequence optimization problem is polynomial time solvable for any $n$,
any graph $H$ on $[n]$ of treewidth bounded by $k$, and any functions $f_1,\dots,f_n$.
\et

We conjecture that degree sequence optimization over the complete graph $H=K_n$ and complete
bipartite graph $H=K_{m,n}$ can be done in polynomial time for all functions $f_i$.

\section{Proof}

We need some notation. Let $H=([n],E)$ be a graph on $[n]$. A {\em nice tree decomposition} of $H$
is a rooted tree $T$ where each node $v\in V(T)$ has a set $I_v\subseteq[n]$ called the {\em bag}
of $v$ such that: every $i\in[n]$ lies in some bag; every edge $\{i,j\}\in E$ is contained in
some bag; if $u,v,w\in V(T)$ where $v$ lies on the path in $T$ between $u$ and $w$, and
$i\in I_u\cap I_w$, then $i\in I_v$; the root $r$ and every leaf $l$ of $T$ satisfy
$I_r=I_l=\emptyset$; every non leaf node $v\in V(T)$ is one of the following:
\begin{itemize}
\item
{\bf introduce node} meaning it has a single child $u$ and $I_v=I_u\uplus\{i\}$ for some $i\in[n]$;
\item
{\bf forget node} meaning it has a single child $u$ and $I_u=I_v\uplus\{i\}$ for some $i\in[n]$;
\item
{\bf join node} meaning it has two children $u,w$ and $I_u=I_v=I_w$.
\end{itemize}
The {\em width} of the tree is the maximum cardinality of a bag minus one. The {\em treewidth} of
$H$, denoted $tw(H)$, is the minimum width of a nice tree decomposition of it. It is known that
for any fixed $k$, if $H$ has treewidth at most $k$, then it is possible to compute in polynomial
time a nice tree decomposition $T$ of $H$ of width at most $k$, with $|V(T)|=O(n)$. See
\cite{CFKLMPPS} for more details. So from here on we may assume that we have a
nice tree decomposition $T$ of $H$.

We let $H[I]$ be the subgraph of $H$ induced by $I\subseteq[n]$.
For $v\in V(T)$ we let $T_v$ be the subtree of $T$ rooted at $v$ and
$I(T_v):=\cup\{I_u:u\in V(T_v)\}$ the union of all bags of nodes in $T_v$.

\vskip.2cm
For every $v\in V(T)$, function $c_v:I_v\rightarrow\{0,1,\dots,n-1\}$,
and subset $F_v\subseteq E(H[I_v])$, let
$$g(v,c_v,F_v)\ :=\ \min\left\{\sum_{i\in I(T_v)}f_i(d_i(G))\ :\
G\subseteq H[I(T_v)],\ E(G[I_v])=F_v,\ d_i(G)=c_v(i),\ i\in I_v\right\}\ ,$$
taken to be $\infty$ if the set is empty. Note that at the tree root $r$ we have
$I(T_r)=[n]$ and so $H[I(T_r)]=H$, and $I_r=\emptyset$ so $F_r$ must be the empty set
and $c_r=c_\emptyset$ the no-value function on $\emptyset$, and so
$g(r,c_\emptyset,\emptyset)=\min\{\sum_{i=1}^nf_i(d_i(G)):G\subseteq H\}$ is
the optimal solution of our problem.

The function $g(v,c_v,F_v)$ is the optimal value of a subproblem of the degree sequence
optimization problem restricted to subgraphs $G$ of the graph $H[I(T_v)]$ induced by the vertex
set $I(T_v)$ which is the union of all bags at nodes of the subtree $T_v$ rooted at $v$,
with further restrictions on $G$ forcing its subset of edges $E(G[I_v])$ induced by the
bag $I_v$ of $v$ to be equal to $F_v$, and the degree $d_i(G)$ of each vertex $i\in I_v$
to be equal to $c_v(i)$. These functions enable us to go down the tree from the leaves to the
root and compute recursively all values $g(v,c_v,F_v)$ at any node $v$ of $T$ by using
all previously computed values $g(u,c_u,F_u)$ at all descendants $u$ of $v$ in the tree.
The proof proceeds by showing, for each of the three types of the non leaf nodes,
how to do these computations. In particular, at any introduce node $v$ or join node $v$,
where these computations are more involved than in a forget node, we will define an auxiliary
function $h(v,c_v,F_v)$, give an explicit formula for it using the previously computed values
$g(u,c_u,F_u)$ at the descendants $u$ of $v$, and then argue that $g(v,c_v,F_v)=h(v,c_v,F_v)$. 
 
We proceed to show how to compute all $g(v,c_v,F_v)$ from the leaves to the root,
thereby solving our degree optimization problem. At any tree leaf $l$ we have
$I(T_l)=I_l=\emptyset$ and $F_l=\emptyset$, $c_l=c_\emptyset$,
and the only graph in the set is $G_\emptyset:=(\emptyset,\emptyset)$
and $g(l,c_\emptyset,\emptyset)=\sum_{i\in\emptyset}f_i(d_i(G_\emptyset))=0$.

\vskip.2cm\noindent
{\bf Introduce node} $v$ with child $u$ and $I_v=I_u\uplus\{i\}$.
For any $c_v$ on $I_v$ and $F_v\subseteq E(H[I_v])$, let $F_u:=F_v\cap E(H[I_u])$, define $c_u$
on $I_u$ by
$$c_u(j)\ :=\ \left\{\begin{array}{ll}
c_v(j), & \hbox{if $\{i,j\}\notin F_v$;} \\
c_v(j)-1, & \hbox{if $\{i,j\}\in F_v$,}
\end{array}\right.$$
and let
$$h(v,c_v,F_v)\ :=\ \left\{\begin{array}{ll}
\infty,\quad
\hbox{if $c_v(i)\neq|\{j\in I_u:\{i,j\}\in F_v\}|$ or $c_u(j)<0$ for some $j\in I_u$;}&\\
g(u,c_u,F_u)+\sum_{j\in I_v}f_j(c_v(j))-\sum_{j\in I_u}f_j(c_u(j)),\quad \hbox{otherwise.} &
\end{array}\right.$$
We now show that $g(v,c_v,F_v)=h(v,c_v,F_v)$. First we show $g(v,c_v,F_v)\geq h(v,c_v,F_v)$.
If $g(v,c_v,F_v)=\infty$ then this clearly holds.
Otherwise pick $G_v\subseteq H[I(T_v)]$ with $E(G_v[I])=F_v$ and $d_j(G_v)=c_v(j)$
for all $j\in I_v$ which attains $\sum_{j\in I(T_v)}f_j(d_j(G_v))=g(v,c_v,F_v)$.
Since $T$ is a tree decomposition of $H$ there is no edge $\{i,j\}$ in $G_v\subseteq H[I(T_v)]$
with $j\notin I_u$ and hence $c_v(i)=|\{j\in I_u:\{i,j\}\in F_v\}|$. Let $G_u=G_v[I(T_u)]$.
Then $E(G[I_u])=F_u$ and for each $j\in I_u$,
$$c_u(j)\ =\ \left\{\begin{array}{ll}
c_v(j)=d_j(G_v)=d_j(G_u), & \hbox{if $\{i,j\}\notin F_v$;} \\
c_v(j)-1=d_j(G_v)-1=d_j(G_u), & \hbox{if $\{i,j\}\in F$,}
\end{array}\right.$$
so $c_u(j)=d_j(G_u)$ and in particular $c_u(j)\geq 0$ for all $j\in I_u$. Therefore
\begin{eqnarray*}
g(v,c_v,F_v)&=&\sum_{j\in I(T_v)}f_j(d_j(G_v))\\
&=&\sum_{j\in I(T_u)}f_j(d_j(G_u))+\sum_{j\in I_v}f_j(c_v(j))-\sum_{j\in I_u}f_j(c_u(j))\\
&\geq&g(u,c_u,F_u)+\sum_{j\in I_v}f_j(c_v(j))-\sum_{j\in I_u}f_j(c_u(j))\\
&=&h(v,c_v,F_v)\ .
\end{eqnarray*}
Next we show $g(v,c_v,F_v)\leq h(v,c_v,F_v)$. If $h(v,c_v,F_v)=\infty$ then we are done. Otherwise
$$h(v,c_v,F_v)\ =\ g(u,c_u,F_u)+\sum_{j\in I_v}f_j(c_v(j))-\sum_{j\in I_u}f_j(c_u(j))$$
and $g(u,c_u,F_u)$ is finite. Pick $G_u\subseteq H[I(T_u)]$ with $E(G_u[I])=F_u$ and
$d_j(G_u)=c_u(j)$ for $j\in I_u$ attaining $\sum_{j\in I(T_u)}f_j(d_j(G_u))=g(u,c_u,F_u)$.
Let $G_v\subseteq H[I(T_v)]$ be such that $G_v[I(T_u)]=G_u$ and $E(G_v)\cap E(H[I_v])=F_v$.
Then $d_j(G_v)=c_v(j)$ for $j\in I_v$, and
\begin{eqnarray*}
g(v,c_v,F_v)&\leq&\sum_{j\in I(T_v)}f_j(d_j(G_v))\\
&=&\sum_{j\in I(T_u)}f_j(d_j(G_u))+\sum_{j\in I_v}f_j(c_v(j))-\sum_{j\in I_u}f_j(c_u(j))\\
&=&g(u,c_u,F_u)+\sum_{j\in I_v}f_j(c_v(j))-\sum_{j\in I_u}f_j(c_u(j))
\ =\ h(v,c_v,F_v)\ .
\end{eqnarray*}
We need to compute $h(v,c_v,F_v)$ for $n^{|I_v|}=O(n^{k+1})$
values of $c_v$ and $O(2^{|I_v|\choose 2})=2^{O(k^2)}$ possible sets $F_v$,
each computation involving a constant number of operations, so all is polynomial.

\vskip.2cm\noindent
{\bf Forget node} $v$ with child $u$ and $I_u=I_v\uplus\{i\}$.
Then for any $c_v$ on $I_v$ and $F_v\subseteq E(H[I_v])$,
$$g(v,c_v,F_v)\ =\ \min\left\{g(u,c_u,F_u):
F_u\subseteq E(H[I_u]),\ F_u\cap E(H[I_v])=F_v,\ c_u(j)=c_v(j),\ j\in I_v\right\}\ .$$
We need to compute this for $n^{|I_v|}=O(n^k)$ values of $c_v$ and
$O(2^{|I_v|\choose 2})=2^{O(k^2)}$ possible $F_v$, each computation involving
taking minimum over a set of size $O(n^2)$, so all is polynomial.

\vskip.2cm\noindent
{\bf Join node} $v$ with children $u,w$ and $I_u=I_v=I_w$.
For any $c_v$ on $I_v$ and $F_v\subseteq E(H[I_v])$ let
\begin{eqnarray*}
h(v,c_v,F_v)
&:=&\min\{g(u,c_u,F_u)+g(w,c_w,F_w)\\
&+&\sum_{i\in I_v}\left(f_i(c_v(i))-f_i(c_u(i))-f_i(c_w(i))\right)
\ :\ F_v=F_u\uplus F_w,\ c_v=c_u+c_w\}\ ,
\end{eqnarray*}
the set running over partitions $F_v=F_u\uplus F_w$ and function decompositions $c_v=c_u+c_w$.

We now show that $g(v,c_v,F_v)=h(v,c_v,F_v)$. Let $I:=I_u=I_v=I_w$.
Let $I^c_u:=I(T_u)\setminus I$ and $I^c_w:=I(T_w)\setminus I$.
First we show $g(v,c_v,F_v)\geq h(v,c_v,F_v)$. If $g(v,c_v,F_v)=\infty$ then this clearly holds.
Otherwise pick $G_v\subseteq H[I(T_v)]$ with $E(G_v[I])=F_v$ and $d_i(G_v)=c_v(i)$
for all $i\in I$ which attains $\sum_{i\in I(T_v)}f_i(d_i(G_v))=g(v,c_v,F_v)$.
Consider any partition $F_v=F_u\uplus F_w$. Let $G_u\subseteq G_v[I(T_u)]$ be the subgraph
with $E(G_u)=E(G_v[I(T_u)])\setminus F_w$ and let $G_w\subseteq G_v[I(T_w)]$ be the subgraph
with $E(G_w)=E(G_v[I(T_w)])\setminus F_u$. Define functions $c_u$ and $c_w$ on $I$ by
$c_u(i):=d_i(G_u)$ and $c_w(i):=d_i(G_w)$ for all $i\in I$.
Since $T$ is a tree decomposition of $H$ we have $I^c_u\cap I^c_w=\emptyset$.
Hence $E(G_v)=E(G_u)\uplus E(G_w)$, and
$d_i(G_v)=d_i(G_u)$ for $i\in I^c_u$, $d_i(G_v)=d_i(G_w)$ for $i\in I^c_w$, and
$d_i(G_v)=d_i(G_u)+d_i(G_w)$ for $i\in I$ so $c_v=c_u+c_w$.
Hence
\begin{eqnarray*}
g(v,c_v,F_v)&=&\sum_{i\in I(T_v)}f_i(d_i(G_v))\\
&=&\sum_{i\in I^c_u}f_i(d_i(G_u))+\sum_{i\in I}f_i(c_v(i))+\sum_{i\in I^c_w}f_i(d_i(G_w))
\ =\ \sum_{i\in I}f_i(c_v(i))\\
&+&\left(\sum_{i\in I(T_u)}f_i(d_i(G_u))-\sum_{i\in I}f_i(c_u(i))\right)
+\left(\sum_{i\in I(T_w)}f_i(d_i(G_w))-\sum_{i\in I}f_i(c_w(i))\right)\\
&=&\sum_{i\in I(T_u)}f_i(d_i(G_u))+\sum_{i\in I(T_w)}f_i(d_i(G_w))
+\sum_{i\in I}\left(f_i(c_v(i))-f_i(c_u(i))-f_i(c_w(i))\right)\\
&\geq&g(u,c_u,F_u)+g(w,c_w,F_w)+\sum_{i\in I_v}\left(f_i(c_v(i))-f_i(c_u(i))-f_i(c_w(i))\right)\\
&\geq&h(v,c_v,F_v)\ .
\end{eqnarray*}
Conversely, consider a partition $F_v=F_u\uplus F_w$ and a decomposition $c_v=c_u+c_w$ attaining
$$g(u,c_u,F_u)+g(w,c_w,F_w)+
\sum_{i\in I_v}\left(f_i(c_v(i))-f_i(c_u(i))-f_i(c_w(i))\right)\ =\ h(v,c_v,F_v)\ .$$
If $h(v,c_v,F_v)=\infty$ then $g(v,c_v,F_v)\leq h(v,c_v,F_v)$.
Otherwise $g(u,c_u,F_u)$ and $g(w,c_w,F_w)$ are finite. Pick
$G_u\subseteq H[I(T_u)]$ with $E(G_u[I])=F_u$ and $d_i(G_u)=c_u(i)$
for $i\in I$ attaining $\sum_{i\in I(T_u)}f_i(d_i(G_u))=g(u,c_u,F_u)$,
and pick $G_w\subseteq H[I(T_w)]$ with $E(G_w[I])=F_u$ and $d_i(G_w)=c_w(i)$
for $i\in I$ attaining $\sum_{i\in I(T_w)}f_i(d_i(G_w))=g(w,c_w,F_w)$.
Let $G_v\subseteq H[I(T_v)]$ be such that $E(G_v)=E(G_u)\uplus E(G_w)$.
Then $E(G_v[I])=F_v$ and $d_i(G_v)=c_v(i)$ for $i\in I$, and
\begin{eqnarray*}
g(v,c_v,F_v)&\leq&\sum_{i\in I(T_v)}f_i(d_i(G_v))\ =\
\sum_{i\in I(T_u)}f_i(d_i(G_u))+\sum_{i\in I(T_w)}f_i(d_i(G_w))\\
&+&\sum_{i\in I}\left(f_i(c_v(i))-f_i(c_u(i))-f_i(c_w(i))\right)\\
&=&g(u,c_u,F_u)+g(w,c_w,F_w)+\sum_{i\in I}\left(f_i(c_v(i))-f_i(c_u(i))-f_i(c_w(i))\right)\\
&=&h(v,c_v,F_v)\ .
\end{eqnarray*}
We need to compute $g(v,c_v,F_v)=h(v,c_v,F_v)$ for $n^{|I_v|}=O(n^{k+1})$ values of $c_v$ and for
$O(2^{|I_v|\choose 2})=2^{O(k^2)}$ possible $F_v$. Each computation involves taking the minimum
over all $2^{O(k^2)}$ partitions $F_v=F_u\uplus F_w$ and $O(n^{k+1})$ decomposition $c_v=c_u+c_w$.
So all is polynomial.

\section*{Acknowledgments}
S. Onn was supported by a grant from the Israel Science Foundation and the Dresner chair.

\end{document}